\documentclass[siam10.clo]{siamltex}

\def\cprime{$'$}

\begin{document} 

\catcode`\Œ=\active\defŒ{{\aa}}       
\catcode`\=\active\def{\c c}        
\catcode`\'=\active\def'{\c C}        

\catcode`\Š=\active\defŠ{\"a}        
\catcode`\'=\active\def'{\"e}        
\catcode`\•=\active\def•{\"{\i}}     
\catcode`\š=\active\defš{\"o}        
\catcode`\Ÿ=\active\defŸ{\"u}        
\catcode`\…=\active\def…{\"O}        
\catcode`\†=\active\def†{\"U}        
\catcode`\‡=\active\def‡{\'a}        
\catcode`\Ž=\active\defŽ{\'e}        
\catcode`\'=\active\def'{\'{\i}}     
\catcode`\—=\active\def—{\'o}        
\catcode`\œ=\active\defœ{\'u}        
\catcode`\ƒ=\active\defƒ{\'E}        
\catcode`\ˆ=\active\defˆ{\`a}        
\catcode`\=\active\def{\`e}        
\catcode`\"=\active\def"{\`{\i}}     
\catcode`\˜=\active\def˜{\`o}        
\catcode`\=\active\def{\`u}        
\catcode`\Ë=\active\defË{\`A}        
\catcode`\‹=\active\def‹{\~a}        
\catcode`\–=\active\def–{\~n}        
\catcode`\›=\active\def›{\~o}        
\catcode`\Ì=\active\defÌ{\~A}        
\catcode`\"=\active\def"{\~N}        
\catcode`\Í=\active\defÍ{\~O}        
\catcode`\‰=\active\def‰{\^a}        
\catcode`\=\active\def{\^e}        
\catcode`\"=\active\def"{\^{\i}}     
\catcode`\™=\active\def™{\^o}        
\catcode`\ž=\active\defž{\^u}        

\font\goth=eufm10 

\newcommand{\Beta}{{\mathrm B}}

\newcommand{\BbbC}{{\mathbf C}}
\newcommand{\BbbF}{{\mathbf F}}
\newcommand{\BbbQ}{{\mathbf Q}}
\newcommand{\BbbR}{{\mathbf R}}
\newcommand{\BbbZ}{{\mathbf Z}}

\newcommand{\Qbar}{\overline{\BbbQ}}

\newcommand{\rmh}{{\mathrm h}}

\def\og{\leavevmode\raise.3ex\hbox{$\scriptscriptstyle
\langle\!\langle$}}

\def\fg{\leavevmode\raise.3ex\hbox{$\scriptscriptstyle
\,\rangle\!\rangle$}}

\newcommand{\calE}{{\mathcal E}}

\def\rmM{{\mathrm M}}

\newcommand{\calP}{{\mathcal P}}

\newcommand{\us}{\underline{s}}

\def\Zgoth{\hbox{\goth Z}}

\newcommand{\calL}{{\mathcal L}}

\newcommand{\bG}{\mathbf{G}}
\newcommand{\bGm}{\bG_m}
\newcommand{\bGa}{\bG_a}

\def\tg#1{T_e(#1)}

\newcommand{\virgule}{\raise 2 pt \hbox{,}}

\newcommand{\calF}{{\mathcal F}}

\newtheorem{theoreme}[equation]{ThŽorme}
\newtheorem{conjecture}[equation]{Conjecture} 
\newtheorem{corollaire}[equation]{Corollaire} 
\newtheorem{lemme}[equation]{Lemme}
\newtheorem{remarque}[equation]{Remarque}
\newtheorem{question}[equation]{Question}
\newtheorem{exemple}[equation]{Exemple}
\newtheorem{propos}[equation]{Proposition}
 
\title{Transcendance de p\'eriodes: \\ \'etat des connaissances
}

\author{Michel Waldschmidt\thanks{Institut de Math\'ematiques 
de Jussieu -- UMR 7586 du CNRS,
Universit\'e P.~et M.~Curie (Paris VI),
175 rue du Chevaleret,
F-75013 Paris
{\tt miw@math.jussieu.fr}
{\tt http://www.math.jussieu.fr/$\sim$miw}}}

\maketitle




\begin{abstract}
Les nombres rŽels ou complexes forment un ensemble ayant la puissance du continu. Parmi eux, ceux qui sont \og intŽressants\fg, qui apparaissent \og naturellement\fg, qui mŽritent notre attention, forment un ensemble dŽnombrable. Dans cet Žtat d'esprit nous nous intŽressons aux pŽriodes au sens de  Kontsevich et Zagier. Nous faisons le point sur l'Žtat de nos connaissances concernant la nature arithmŽtique de ces nombres: dŽcider si une pŽriode est un nombre rationnel, algŽbrique irrationnel ou au contraire transcendant est  l'objet de quelques thŽormes et de beaucoup de conjectures. Nous prŽcisons aussi ce qui est connu sur l'approximation diophantienne de tels nombres, par des nombres rationnels ou algŽbriques.

\end{abstract}

\pagestyle{myheadings}
\thispagestyle{plain}
\markboth{Michel Waldschmidt}{\small Transcendance de p\'eriodes
}


\section{Introduction}

Dans leur article  \cite{KontsevichZagierMU-00} intitulŽ \og Periods\fg, M.~Kontsevich et D.~Zagier introduisent la notion de pŽriodes en en donnant deux dŽfinitions dont ils disent qu'elles sont Žquivalentes; il proposent une conjecture, deux principes et cinq problmes. Le premier principe est le suivant: \og {\it chaque fois que vous rencontrez un nouveau nombre et que vous voulez savoir s'il est transcendant, commencez par essayer de savoir si c'est une pŽriode}\fg. 

Si la rŽponse est nŽgative, alors le nombre est transcendant; en effet les pŽriodes forment une sous algbre de $\BbbC$ sur le corps $\Qbar$ des nombres algŽbriques, donc tout nombre algŽbrique est une pŽriode. 

Le but de cet exposŽ est d'examiner ce qui se passe si la rŽponse est positive: {\it que sait-on sur la transcendance de pŽriodes? }

Nous considŽrons aussi l'aspect quantitatif de cette question, en liaison avec la question suivante de  \cite{KontsevichZagierMU-00}, \S~1.2  qui prŽcde leur conjecture 1: quand on veut vŽrifier une ŽgalitŽ entre deux nombres algŽbriques, il suffit de calculer ces deux nombres avec une prŽcision suffisante, puis d'utiliser l'inŽgalitŽ de Liouville qui Žtablit que deux nombres algŽbriques distincts de degrŽ et hauteur bornŽe ne peuvent tre trop proches l'un de l'autre. Dans l'exemple qu'ils donnent, dž ˆ D.~Shanks \cite{ShanksFQ-74}:
\begin{equation}
\label{E:Shanks}
\sqrt{11+2\sqrt{29}}+
\sqrt{16-2\sqrt{29}+2\sqrt{55-10\sqrt{29}}}=
\sqrt{5}+\sqrt{22+2\sqrt{5}},
\end{equation}
la diffŽrence $\gamma$ entre les deux membres de (\ref{E:Shanks}) est un nombre algŽbrique de degrŽ $\le 16$ sur $\BbbQ(\sqrt{5},\sqrt{29})$, donc de degrŽ $\le 64$ sur $\BbbQ$. Pour chacun des  $64$ ŽlŽments $\epsilon=(\epsilon_1,\ldots,\epsilon_6)\in\{0,1\}^6$, posons
\begin{eqnarray}
\gamma_\epsilon=
\epsilon_3\sqrt{11+2\epsilon_2\sqrt{29}}
&+&
\epsilon_4\sqrt{16-2\epsilon_2\sqrt{29}+2\epsilon_5\sqrt{55-10\epsilon_2\sqrt{29}}}\nonumber\\
&+&\epsilon_1\sqrt{5}+\epsilon_6\sqrt{22+2\epsilon_1\sqrt{5}}.
\nonumber
\end{eqnarray}
Le nombre
$$
N= \prod_{\epsilon}\gamma_\epsilon
$$
est un entier rationnel. Il suffit de le calculer avec une prŽcision d'un chiffre aprs la virgule pour vŽrifier qu'il satisfait $-1<N<1$, donc qu'il est nul (c'est le cas le plus simple de l'inŽgalitŽ de Liouville  \cite{miwGL326-00} \S~3.5). Il s'ensuit qu'un (au moins) des $64$ facteurs $\gamma_\epsilon$ du produit est nul, et l'ŽgalitŽ (\ref{E:Shanks}) s'en dŽduit aisŽment.

 La question posŽe par Kontsevich et Zagier dans \cite{KontsevichZagierMU-00} \S~1.2 consiste ˆ savoir si on peut faire de mme avec les pŽriodes. Il s'agirait de dŽfinir une notion  de complexitŽ d'une pŽriode analogue ˆ celle de hauteur pour un nombre algŽbrique, puis de minorer cette complexitŽ pour une pŽriode non nulle afin de remplacer l'inŽgalitŽ de Liouville.  Une des suggestions qu'ils font est de compter le nombre de touches nŽcessaires pour taper en \TeX\ une intŽgrale dont la valeur est la pŽriode en question. 
 
 Dans cet Žtat d'esprit il serait intŽressant de savoir s'il existe des nombres qui sont ˆ la fois une pŽriode et un nombre de Liouville. Une rŽponse nŽgative signifierait que {\it pour toute pŽriode rŽelle $\theta$, il existe une constante $c(\theta)>0$ telle que, pour tout nombre rationnel $p/q$ distinct de $\theta$ avec $q\ge 2$, on ait}
 $$
 \left|\theta-\frac{p}{q}\right|>\frac{1}{q^{c(\theta)}}\cdotp
 $$
 Plus ambitieusement on peut demander si les pŽriodes (complexes) se comportent, pour l'approximation par des nombres algŽbriques, comme presque tous les nombres (complexes) \cite{BugeaudCTM-2004, miwGL326-00}: {\it Žtant donnŽe une pŽriode transcendante $\theta\in\BbbC$, existe-t-il une constante $\kappa(\theta)$ telle que, pour tout polyn™me non nul $P\in\BbbZ[X]$, on ait
 $$
 |P(\theta)|\ge H^{-\kappa(\theta) d},
 $$
 o $H\ge 2$ est un majorant de la hauteur (usuelle) de $P$ (maximum des valeurs absolues des coefficients) et $d$ son degrŽ?}

\section{IntŽgrales abŽliennes}

La nature arithmŽtique de la valeur de l'intŽ\-grale d'une fonction algŽbrique d'une variable entre des bornes algŽbriques (ou infinies) est maintenant bien connue, aussi bien sous l'aspect qualitatif que quantitatif.


\subsection{Genre $0$: logarithmes de nombres algŽbriques}

L'outil principal est le thŽorme de Baker sur l'indŽpendance linŽaire, sur le corps $\Qbar$ des nombres algŽbriques, de logarithmes de nombres algŽbriques. Nous n'utilisons ici que le cas particulier suivant:

\begin{theoreme}
\label{T:Baker}
Soient $\alpha_1,\ldots,\alpha_n$ des nombres algŽbriques non nuls, $\beta_1,\ldots,\beta_n$ des nombres algŽbriques, et, pour $1\le i\le n$, $\log\alpha_i$ un logarithme complexe de $\alpha_i$. Alors le nombre 
$$
\beta_1\log\alpha_1+\cdots+\beta_n\log\alpha_n
$$
est soit nul, soit transcendant.

\end{theoreme}

On en dŽduit:

\begin{corollaire}
\label{C:genrezero}
Soient $P$ et $Q$ des polyn™mes ˆ coefficients algŽbriques vŽrifiant  
$\deg P<\deg Q$ et soit $\gamma$ un chemin fermŽ, ou bien un chemin dont les extrŽmitŽs sont algŽbriques ou infinies. Si l'intŽgrale 
\begin{equation}\label{E:integraledePsurQ}
\int_\gamma \frac {P(z)}{Q(z)} dz
\end{equation}
existe, alors elle est soit nulle, soit transcendante.
\end{corollaire}

 Un exemple cŽlbre \cite{SiegelAMS-49} p.~97 est 
 $$
\int_0^1\frac{dt}{1+t^3}=\frac{1 }{
3}\left(\log2+\frac{\pi}{\sqrt{3}}\right)
\cdotp
$$

Le corollaire \ref{C:genrezero} se dŽduit du thŽorme \ref{T:Baker} en dŽcomposant la fraction rationnelle $P(z)/Q(z)$ en ŽlŽments simples
 (voir  par exemple \cite{vanderPoortenPAMS-71}). En fait le corollaire  \ref{C:genrezero} est Žquivalent au thŽorme \ref{T:Baker}: il suffit d'Žcrire le logarithme d'un nombre algŽbrique comme une pŽriode;  pour la dŽtermination principale, quand $\alpha$ n'est pas rŽel nŽgatif, on a par exemple
 $$
 \log\alpha=\int_0^\infty \frac{(\alpha-1)dt}{(t+1)(\alpha t+1)}
  $$
tandis que
 $$
 i\pi=2i\int_0^\infty \frac{dt}{1+t^2}\cdotp
 $$ 
 Les mesures d'indŽpendance linŽaire de logarithmes de nombres algŽbriques (minorations de combinaisons linŽaires, ˆ coefficients algŽbriques, de lo\-ga\-ri\-thmes de nombres algŽbriques - voir par exemple \cite{miwGL326-00}) contiennent le fait qu'une intŽgrale non nulle de la forme (\ref{E:integraledePsurQ}) a une valeur absolue minorŽe explicitement en termes des hauteurs de $P$ et $Q$ et de leurs degrŽs, ainsi que des hauteurs et degrŽs des nombres algŽbriques extrŽmitŽs de $\gamma$.


\subsection{Genre $1$: intŽgrales elliptiques}

La nature arithmŽtique des va\-leurs d'intŽgrales elliptiques de premire ou deuxime espce a ŽtŽ ŽtudiŽe ds 1934 \cite{SchneiderCrelle-34} puis 1937 \cite{SchneiderMA-37} par Th.~Schneider. Voici le thŽorme 15 version III de \cite{Schneider-57}.

\begin{theoreme}\label{T:Schneider-genre-un}
Toute intŽgrale elliptique de premire ou deuxime espce ˆ coefficients algŽbriques et calculŽe entre des bornes algŽbriques distinctes a pour valeur un nombre nul ou transcendant.

\end{theoreme}

En particulier {\it toute pŽriode non nulle d'une intŽgrale elliptique de premire ou deuxime espce ˆ coefficients algŽbriques est transcendante.}

Le  thŽorme 16   de \cite{Schneider-57} concerne  {\it  la transcendance du quotient de deux intŽgrales elliptiques de premire espce. }

Une consŽquence que cite Schneider de son thŽorme 17 dans  \cite{Schneider-57}  s'Žnonce:   {\it  la valeur prise par une intŽgrale elliptique de premire ou de deuxime espce ˆ coefficients algŽbriques entre des bornes algŽbriques est quotient d'une pŽriode par un facteur rationnel ou transcendant. }

Du thŽorme \ref{T:Schneider-genre-un} on dŽduit le rŽsultat citŽ dans \cite{KontsevichZagierMU-00} \S~1.1:   {\it   si $a$ et  $b$ sont deux nombres algŽbriques rŽels positifs, l'ellipse  dont les longueurs d'axes sont $a$ et $b$ a un pŽrimtre 
\begin{equation}
\label{E:perimetre}
2\int_{-b}^{b} \sqrt{1+\frac{a^2 x^2}{ b^4-b^2x^2}}\ dx
\end{equation}
qui est  un nombre transcendant. }  Plus gŽnŽralement   {\it   la longueur de tout arc dont les extrŽmitŽs sont des points de coordonnŽes algŽbriques est un nombre transcendant ou nul.}

 Il en est de mme pour une lemniscate 
 $$
(x^2+ y^2)^2=2a^2(x^2-y^2)
$$
quand $a$ est algŽbrique.

 Ces ŽnoncŽs sont dŽmontrŽs par Schneider comme consŽquences de rŽsultats sur les fonctions elliptiques.  Voici par exemple la version I du  thŽorme 15 de \cite{Schneider-57}. Soit $\wp$ une fonction elliptique de Weierstrass d'invariants $g_2$ et $g_3$ algŽbriques:
 $$
 {\wp'}^2=4\wp^3-g_2\wp-g_3.
 $$
 Soient $\zeta$ la fonction zta de Weierstrass associŽe ˆ $\wp$,  $a$ et $b$ deux nombres algŽbriques non tous deux nuls et  $u$ un nombre complexe non p™le de $\wp$. 
  {\it  Alors l'un au moins des deux nombres $\wp(u)$, $au+b\zeta(u)$ est transcendant.}

Ainsi en considŽrant les deux courbes elliptiques
$$
y^2=x^3-x \qquad \hbox{et}\qquad
y^2=x^3-x
$$
on en dŽduit que   {\it  chacun des deux nombres
 \begin{equation}
 \label{E:GammaUnQuart}
\int_0^1 \frac{dx }{ \sqrt{x-x^3} }
=\frac{1 }{ 2} \Beta(1/4,1/2)=\frac{
\Gamma(1/4)^2 }{ 2^{3/2}\pi^{1/2}}
\end{equation}
et
 \begin{equation}
 \label{E:GammaUnTiers}
\int_0^1 \frac{dx }{ \sqrt{1-x^3} }
=\frac{1 }{ 3} \Beta(1/3,1/2)=\frac{
\Gamma(1/3)^3 }{  2^{4/3}3^{1/2} \pi}
\end{equation}
est transcendant.}

Ces deux formules (comparer avec \cite{Masser} p.21) sont des cas particuliers de la   formule de Chowla-Selberg  (cf 
 \cite{GrossRohrlichIM-78} et \cite{KontsevichZagierMU-00} \S~2.3) qui exprime les pŽriodes de courbes elliptiques de type CM comme des produits de valeurs de la fonction Gamma d'Euler dont une des dŽfinitions est:
\begin{equation}
\label{E:GammaProduit}
\Gamma(z)=e^{-\gamma z} z^{-1} \prod_{n=1}^\infty \left(
1+\frac{z}{n}\right)^{-1}e^{z/n}.
\end{equation}
L'extension par G.~Shimura aux variŽtŽs abŽliennes de type CM de la formule de Chowla-Selberg donne lieu aux relations de Deligne-Koblitz-Ogus sur la fonction Gamma (voir \cite{BrownawellPapanikolasCrelle-02}).

 
\subsection{Genre $\ge 1$: intŽgrales abŽliennes}

Dans \cite{SchneiderCrelle-41}, Th.~Schneider  Žtend ses rŽsultats aux intŽgrales abŽliennes. La dŽmonstration est une extension en plusieurs variables de ses rŽsultats antŽrieurs; dans la partie analytique de la dŽmonstration de transcendance, l'outil essentiel, un lemme de Schwarz, est Žtendu en plusieurs variables gr‰ce ˆ une formule d'interpolation pour les produits cartŽsiens. Cela permet ˆ Schneider d'obtenir des ŽnoncŽs sur les fonctions abŽliennes. L'exemple le plus important des rŽsultats qu'il obtient est le suivant:

\begin{theoreme}\label{T:Beta}
Soient
 $a$ et  $b$ des nombres rationnels non entiers tels que  $a+b$ ne soit pas non plus un entier. Alors le nombre
\begin{equation}
\label{E:Beta}
\Beta(a,b)=\frac{\Gamma(a)\Gamma(b) }{ \Gamma(a+b)}
=\int_0^1 x^{a-1}(1-x)^{b-1}dx
\end{equation}
est transcendant.
\end{theoreme}

 Les travaux sur la nature arithmŽtique des valeurs d'intŽgrales abŽliennes sont nombreux: ceux de Th.~Schneider en 1940 ont ŽtŽ poursuivis par S.~Lang dans les annŽes 1960, puis notamment par  D.W.~Masser gr‰ce ˆ la mŽthode de Baker dans les annŽes 1980, pour arriver ˆ une solution essentiellement complte de la question en 1989 par 
 G.~WŸstholz \cite{WustholzAnnM-89} qui a obtenu une extension satisfaisante du 
 thŽorme  \ref{T:Baker} de  Baker  aux groupes algŽbriques commutatifs. On conna"t donc essentiellement ce que l'on souhaite sur la  transcendance et  l'indŽpendance linŽaire (sur le corps des nombres algŽbriques) d'intŽgrales abŽliennes de premire, seconde ou troisime espce. Par exemple J.~Wolfart et G.~WŸstholz \cite{WolfartWustholzMA-85} ont montrŽ que les seules relations linŽaires ˆ coefficients algŽbriques entre les valeurs $\Beta(a,b)$ de la fonction Bta en des points $(a,b)\in\BbbQ^2$ sont celles qui rŽsultent des relations de Deligne-Koblitz-Ogus.
 
 De plus on dispose Žga\-lement maintenant de rŽsultats quantitatifs qui permettent de minorer la valeur d'une intŽgrale abŽlienne quand elle est non nulle -- les estimations les plus rŽcentes et les plus prŽcises sur ce sujet, dans le cadre gŽnŽral des groupes algŽbriques, sont dues ˆ ƒ.~Gaudron \cite{GaudronCRAS-01, GaudronFLLEVA-04}.

 Si les relations linŽaires ˆ coefficients algŽbriques entre les valeurs d'intŽ\-grales abŽliennes sont maintenant bien connues, il n'en est pas de mme des relations algŽbriques. Dans une note de bas de page \cite{GrothendieckIHES-66}, A.~Grothendieck propose un ŽnoncŽ conjectural sur la transcendance de pŽriodes de variŽtŽs abŽliennes dŽfinies sur le corps des nombres algŽbriques. La premire formulation prŽcise de cette conjecture est donnŽe par S.~Lang dans son livre \cite{LangITN-66}, o l'on trouve aussi la premire formulation de la conjecture de Schanuel sur l'indŽpendance algŽbrique des valeurs de la fonction exponentielle (voir aussi \cite{DeligneMSMF-80}). Ces ŽnoncŽs ont ŽtŽ dŽveloppŽs par Y.~AndrŽ (voir notamment \cite{AndreQCTIGA-97}) qui propose une gŽnŽralisation commune des conjectures de Grothendieck et Schanuel. Pour des $1$-motifs attachŽs aux produits de courbes elliptiques,  C.~Bertolin \cite{BertolinJNT-02} a explicitŽ la situation conjecturale en formulant sa conjecture elliptico-torique qui fait intervenir la fonction exponentielle, les fonctions $\wp$ et $\zeta$ de Weierstrass, les intŽgrales elliptiques et l'invariant modulaire $j$   (voir aussi \cite{miwODPJMS-04} et \cite{miwHyderabad-03}).

\section{Valeurs de la fonction Gamma d'Euler}
La dŽfinition (\ref{E:Beta}) de la fonction Bta  sous forme d'une intŽgrale montre que ses valeurs aux points de $\BbbQ^2$ o elle est dŽfinie sont des pŽriodes. De la relation (\ref{E:Beta}) entre les fonctions Gamma et Bta on dŽduit
$$
\Gamma(a_1)\cdots \Gamma(a_n) = \Gamma(a+\cdots+a_n)
\prod_{i=1}^{n-1} \Beta(a_1+\cdots+a_{i-1},a_i).
$$
Il en rŽsulte que pour tout $p/q\in\BbbQ$ avec $p>0$ et $q>0$, le nombre $\Gamma(p/q)^q$ est une pŽriode. Par exemple
$$
\pi=\Gamma(1/2)^2=\int_0^1x^{-1/2}(1-x)^{-1/2}dx.
$$
De  (\ref{E:GammaUnQuart}) et (\ref{E:GammaUnTiers}) on dŽduit aussi des expressions de $\Gamma(1/3)^3$ et $\Gamma(1/4)^4$ comme pŽriodes.

On conna"t bien mieux la nature arithmŽtique des valeurs de la fonction  Bta  d'Euler (gr‰ce au thŽorme \ref{T:Beta} de Schneider) que celles de la fonction Gamma. On sait que le nombre $\Gamma(1/2)=\sqrt{\pi}$ est transcendant, gr‰ce ˆ Lindemann. La transcendance de $\Gamma(1/4)$ et $\Gamma(1/3)$ a ŽtŽ Žtablie par  G.V.~{\v{C}}udnovs{\cprime}ki\u{\i}  \cite{ChudnovskyDANU-76}.
   
\begin{theoreme} \label{T:Chudnovsky}
Les deux nombres
$$
\Gamma(1/4)\quad\hbox{ et}\quad \pi
$$
sont algŽbriquement indŽpendants, et il en est de mme des deux nombres
$$
\Gamma(1/3)\quad\hbox{et}\quad \pi.
$$

\end{theoreme}

Comme l'a remarquŽ D.W. Masser on peut aussi Žnoncer ces rŽsultats en disant que les deux nombres
$\Gamma(1/4)$ et $\Gamma(1/2)$
sont algŽbriquement indŽpendants et qu'il en est de mme des deux nombres
$\Gamma(1/3)$ et $ \Gamma(2/3)$.

 Les seules autres valeurs de la fonction $\Gamma$ en des points rationnels dont on sache dŽmontrer la transcendance sont celles que l'on dŽduit de la transcendance en $1/2$, $1/3$ et $1/4$ en utilisant les relations standard satisfaites par la fonction Gamma (voir ci-dessous). Par exemple $\Gamma(1/6)$ est aussi un nombre transcendant.

La dŽmonstration par G.V.~{\v{C}}udnovs{\cprime}ki\u{\i}
de son thŽorme \ref{T:Chudnovsky} repose sur le  rŽsultat suivant \cite{ChudnovskyDANU-76} concernant les pŽriodes et quasi pŽriodes de fonctions de Weierstrass, que l'on applique aux courbes elliptiques $y^2=x^3-x$ et $y^2=x^3-1$ gr‰ce ˆ (\ref{E:GammaUnQuart}) et (\ref{E:GammaUnTiers})

 {\it 
Soit $\wp$ une fonction elliptique de Weierstrass d'invariants
 $g_2$ et  $g_3$. Soit $\omega$  une pŽriode non nulle de $\wp$ et soit $\eta$  la quasi-pŽriode associŽe de la fonction zta de Weierstrass
 $\zeta$:
 $$
{\wp'}^2=4\wp^3-g_2\wp-g_3,
\qquad
\zeta'=-\wp,\quad
\zeta(z+\omega)=\zeta(z)+\eta.
$$
Alors deux au moins des nombres
$$
g_2,\
g_3,\
\omega/\pi,\
\eta/\pi
$$
sont algŽbriquement indŽpendants.}

De cet ŽnoncŽ on dŽduit aussi le fait que le nombre (\ref{E:perimetre}) est non seulement transcendant, mais mme algŽbriquement indŽpendant de $\pi$ (cf. \cite{KontsevichZagierMU-00} \S~1.1).

Pour l'instant on ne sait pas obtenir la transcendance (sur $\BbbQ$) de $\Gamma(1/4)$ ni de $\Gamma(1/3)$ sans Žtablir le rŽsultat plus fort qui est la transcendance de chacun de ces nombres  sur le corps $\BbbQ(\pi)$. D'un point de vue quantitatif de bonnes mesures de transcendance de ces nombres ont ŽtŽ Žtablies par P.~Philippon puis S.~Bruiltet \cite{BruiltetAA-02}:

\begin{theoreme}
\label{T:Bruiltet}
Pour un polyn™me non constant $P\in\BbbZ[X,Y]$ de degrŽ $d$ et de hauteur $H$, on a
$$
\log|P(\pi,\Gamma(1/4)|>-10^{326} \bigl((\log H+d\log(d+1)\bigr)
d^2\bigl(\log (d+1)\bigr)^2
$$
et
$$
\log|P(\pi,\Gamma(1/3)|>-10^{330} \bigl((\log H+d\log(d+1)\bigr)
d^2\bigl(\log (d+1)\bigr)^2.
$$

\end{theoreme}

Ainsi $\Gamma(1/4)$ et $\Gamma(1/3)$ ne sont pas des nombres de Liouville.

La prochaine Žtape pourrait tre la transcendance du nombre $\Gamma(1/5)$ (cf. \cite{Masser}, p.~2 et p.~35). Du thŽorme \ref{T:Beta} de Schneider sur la transcendance du nombres $\Beta(1/5,1/5)$ on dŽduit que l'un au moins des deux nombres $\Gamma(1/5)$, $\Gamma(2/5)$ est transcendant. Un rŽsultat plus prŽcis se dŽduit des travaux de P.~Grinspan \cite{GrinspanJNT-02} (voir aussi \cite{VasilevMZ-96}):

\begin{theoreme}
\label{T:Grinspan}
Un au moins des deux  nombres $\Gamma(1/5)$, $\Gamma(2/5)$  est transcendant sur le corps $\BbbQ(\pi)$.

\end{theoreme}

Autrement dit, deux au moins des trois nombres $\Gamma(1/5)$, $\Gamma(2/5)$ et $\pi$ sont algŽbriquement indŽpendants.  La dŽmonstration de  \cite{GrinspanJNT-02} fournit de plus un rŽsultat quantitatif. 

Comme la courbe de Fermat 
$x^5+y^5=z^5$ d'exposant $5$ est de genre $2$, sa jacobienne est une surface abŽlienne; il faut donc remplacer dans la dŽmonstration de {\v{C}}udnovs{\cprime}ki\u{\i}  les fonctions elliptiques par des fonctions abŽliennes, et c'est pourquoi il est difficile de sŽparer les deux nombres $\Gamma(1/5)$ et $\Gamma(2/5)$ quand on veut obtenir la transcendance de chacun d'eux.

Avant de poursuivre avec le dŽnominateur $5$, revenons aux dŽnominateurs $3$ et $4$. Le thŽorme \ref{T:Chudnovsky} 
a ŽtŽ Žtendu par Yu.V.~Nesterenko  \cite{NesterenkoMS-96, lnm1752}, qui obtient l'indŽpendance algŽbrique de trois nombres: 

\begin{theoreme}
\label{T:Nesterenko}   
Les trois nombres
$$
\Gamma(1/4),\quad\pi \quad\hbox{et}\quad e^\pi
$$
sont algŽbriquement indŽpendants, et il en est de mme de
$$
\Gamma(1/3),\quad\pi \quad\hbox{et}\quad e^{\pi\sqrt{3}}.
$$
\end{theoreme}

La dŽmonstration par Yu.V.~Nesterenko de son thŽorme \ref{T:Nesterenko} utilise les sŽries d'Eisenstein $E_2$, $E_4$ et $E_6$ (nous utilisons les notations $P$, $Q$, $R$ de Ramanujan):

\begin{equation}
\label{E:Ramanujan}
\left\{
\begin{array}{cccccc}
P(q)\displaystyle
=E_2(q)=1-24\sum_{n=1}^\infty\frac{nq^n }{ 1-q^n}\virgule\\
Q(q)\displaystyle
=E_4(q)=1+240\sum_{n=1}^\infty\frac{n^3q^n }{ 1-q^n}\virgule\\
R(q)\displaystyle
=E_6(q)=1-504\sum_{n=1}^\infty\frac{n^5q^n }{ 1-q^n}\cdotp\\
\end{array} \right.
\end{equation}

Les premiers rŽsultats de transcendance sur les valeurs de ces fonctions sont dus ˆ D.~Bertrand dans les annŽes 70. Une avancŽe remarquable a ŽtŽ faite  en 1996 par K.~BarrŽ-Sirieix, G.~Diaz, F.~Gramain et G.~Philibert \cite{BDGP-IM-96}, qui ont rŽsolu le problme suivant de Manin (et de Mahler dans le cas $p$-adique) concernant la fonction modulaire $J=Q^3/\Delta$, o $\Delta=12^{-3}(Q^3-R^2)$: {\it pour tout $q\in\BbbC$ avec $0<|q|<1$, l'un au moins des deux nombres $q$, $J(q)$ est transcendant.}

C'est cette percŽe qui a permis ˆ Yu.V.~Nesterenko \cite{NesterenkoMS-96} de dŽmontrer le rŽsultat suivant, citŽ dans le \S~2.4 de \cite{KontsevichZagierMU-00}:

 {\it Soit $q\in\BbbC$  un nombre complexe satisfaisant
$0<|q|<1$. Alors trois au moins des quatre nombres
$$ 
q,\
 P(q),\
  Q(q),\
   R(q)
$$
sont algŽbriquement indŽpendants.
}

Le thŽorme \ref{T:Nesterenko} en rŽsulte en spŽcialisant $q=e^{-2\pi}$ et $q=-e^{-\pi\sqrt{3}}$ car $J(e^{-2\pi})=1728$,
 $$
P(e^{-2\pi})=\frac{3 }{\pi}\virgule
\quad
Q(e^{-2\pi})=3\left(\frac{\omega }{\pi}\right)^4,
\quad
R(e^{-2\pi})=0
$$
avec   (cf. (\ref{E:GammaUnQuart}))
$$
\omega=\frac{\Gamma(1/4)^2 }{\sqrt{8\pi}}=2.6220575542\dots
$$  
tandis que  $J(-e^{-\pi\sqrt{3}})=0$, 
$$
P(-e^{-\pi\sqrt{3}})=\frac{2\sqrt{3} }{\pi}\virgule
\quad
Q(-e^{-\pi\sqrt{3}})=0,
\quad
R(-e^{-\pi\sqrt{3}})=\frac{27 }{
2}\left(\frac{\omega' }{\pi}\right)^6
$$ 
avec (cf. (\ref{E:GammaUnTiers}))
 $$
\omega'=\frac{\Gamma(1/3)^3 }{ 2^{4/3}\pi}=2.428650648\dots
$$  
On trouve dans \cite{KontsevichZagierMU-00} \S~2.3 des commentaires sur les liens entre les pŽriodes et les sŽries d'Eisenstein  (et aussi les fonctions thta, qui interviennent Žga\-lement dans le travail  \cite{NesterenkoMS-96} de Nesterenko --- voir  \cite{lnm1752}).

Le thŽorme \ref{T:Nesterenko} de Nesterenko et celui  \ref{T:Grinspan} de Grinspan 
suggrent le pro\-blme ouvert suivant:

  \begin{conjecture}
  \label{C:GammaDeuxCinquieme}
Trois au moins des quatre nombres
$$
\Gamma(1/5),\quad\Gamma(2/5),\quad\pi \quad\hbox{et}\quad e^{\pi\sqrt{5}}
$$
sont algŽbriquement indŽpendants.
  \end{conjecture}
  
Ce problme fait l'objet de travaux rŽcents de F.~Pellarin (voir en particulier \cite{PellarinIdeauxStables-04}). 

Plus ambitieusement on peut demander quelles sont toutes les relations algŽbriques liant les valeurs de la fonction Gamma en des points rationnels. La question des relations multiplicatives a ŽtŽ considŽrŽe par D.~Rohrlich. On rappelle dŽjˆ ce que sont les relations standard: pour $a\in\BbbC$ (en dehors des p™les de $\Gamma(x)$, $\Gamma(x+1)$, $\Gamma(1-x)$ ou $\Gamma(nx)$ pour que les formules aient un sens), on a
$$
\Gamma(a+1)=a\Gamma(a), \leqno{\hbox{\rm (Translation)}}
$$
$$
\Gamma(a)\Gamma(1-a)=\frac{\pi }{ \sin(\pi a)} \leqno\hbox{{\rm (Reflexion)}}
$$
et, pour tout $n$ entier positif,
$$
\prod_{k=0}^{n-1} \Gamma\left(a+\frac{k }{ n}\right)
=(2\pi)^{(n-1)/2} 
n^{-na+(1/2)}\Gamma(na).
\leqno\hbox{{\rm (Multiplication)  }}
$$

Voici la conjecture de Rohrlich:

\begin{conjecture}
\label{C:Rohrlich}
Toute relation multiplicative de la forme
 $$
 \pi^{b/2}\prod_{a\in\BbbQ}\Gamma(a)^{m_a}\in \Qbar
 $$
 avec $b$ et $m_a$  dans  $\BbbZ$ se dŽduit des relations standard.
 \end{conjecture}

Une formalisation de cette conjecture utilisant la notion de \og distribution universelle\fg\ est donnŽe par S.~Lang dans \cite{LangDPP-77}.

Une conjecture plus ambitieuse que \ref{C:Rohrlich}  est celle de Rohrlich-Lang qui concerne non seulement les relations monomiales, mais plus gŽnŽralement les relations polynomiales: elle prŽtend que  {\it   l'idŽal sur  $\Qbar$  de toutes les relations algŽbriques entre les valeurs de  $(1/\sqrt{2\pi})\Gamma(a)$ pour  $a\in\BbbQ$ est engendrŽ par les relations de distributions, l'Žquation fonctionnelle et l'imparitŽ}.

\section{SŽries de fractions rationnelles} 

Soient $P$ et $Q$ deux fractions rationnelles ˆ coefficients rationnels avec $\deg Q\ge \deg P+2$. {\it Quelle est la nature arithmŽtique de la somme de la sŽrie
\begin{equation}
\label{E:serie}
  \sum_{  {n\ge 0}\atop {Q(n)\not=0}}
  \frac{P(n) }{ Q(n)} \ \hbox{?}
\end{equation}
}

Cette question a ŽtŽ ŽtudiŽe notamment dans
\cite{AdhikariSaradhaShoreyTijdemanIM-01}.
 
La somme de la sŽrie (\ref{E:serie}) peut tre rationnelle: c'est le cas des sŽries {\it tŽlŽscopiques} dont voici des exemples.

\begin{lemme}
\label{L:telescopique}
Soient
$a$ et $b$ deux ŽlŽments de $\BbbC^\times$ tels que $b/a\not\in\BbbZ_{\le 0}$ et soit $k$ un entier $\ge 2$. Alors
\begin{equation}
\label{E:telescopique}
\sum_{n=0}^\infty
\prod_{j=0}^{k-1}
\frac{1}{(an+b+ja)} 
=\frac{1}{(k-1)a}
\prod_{i=0}^{k-2}
\frac{1}{ia+b}\cdotp
\end{equation}

\end{lemme}

En particulier, sous les hypothses du lemme \ref{L:telescopique}, si $a$ et $b$ sont des nombres rationnels alors la sŽrie a pour valeur un nombre rationnel. Ainsi
$$
\sum_{n=2}^{\infty} 
\sum_{m=2}^{\infty} \frac{1}{n^m}=
\sum_{n=2}^{\infty} \frac{1}{n(n-1)}=1.
$$
Un autre  exemple est la somme de la sŽrie (\ref{E:serie}) avec $P(X)=1$ et $Q(X)=(X+1)\cdots(X+k)$ pour $k\ge 2$, ˆ savoir
$$
\sum_{n=0}^\infty
\frac{n!}{(n+k)!} 
=\frac{1}{k-1} \cdot \frac{1}{(k-1)!} \cdotp
$$ 

\noindent{\it DŽmonstration du lemme \ref{L:telescopique}}. 
On utilise la dŽcomposition en ŽlŽments simples de la fraction rationnelle
\begin{equation}
\label{E:ElementsSimples}
\prod_{j=0}^{h}
\frac{1}{X+ja}=\frac{1}{a^h}\sum_{i=0}^h \frac{(-1)^i}{i!(h-i)!}\cdot
\frac{1}{X+ia}
\end{equation}
d'abord avec $h=k-1$: la somme
$$
S=
\sum_{n=0}^\infty
\prod_{j=0}^{k-1}
\frac{1}{(an+b+ja)} 
$$
de la sŽrie du membre de gauche de (\ref{E:telescopique}) s'Žcrit
$$
S=
\sum_{n=0}^\infty
\sum_{i=0}^{k-1} \frac{c_i}{an+b+ia}\ 
$$
avec
$$
c_i=\frac{(-1)^i}{a^{k-1}i!(k-1-i)!}
\qquad
(0\le i\le k-1).
$$
Comme
$$
\sum_{i=0}^{k-1} c_i=0,
$$
on en dŽduit
$$
S=\frac{1}{a^{k-1}}
\sum_{m=0}^{k-2}\frac{d_m}{am+b}
$$
avec
$$
d_m=
\sum_{i=0}^m  \frac{(-1)^i}{i!(k-1-i)!}\cdotp
$$
On vŽrifie par rŽcurrence sur $m$, pour $0\le m\le k-2$,
$$
d_m=\frac{1}{k-1}\cdotp \frac{(-1)^m}{m!(k-2-m)!}\cdotp
$$
Il ne reste plus qu'ˆ appliquer de nouveau (\ref{E:ElementsSimples}), mais cette fois-ci avec $h=k-2$.

 \bigskip

La somme d'une sŽrie (\ref{E:serie}) peut aussi tre transcendante: des exemples 
\cite{AdhikariSaradhaShoreyTijdemanIM-01} 
sont
$$
 \sum_{n=0}^\infty\frac{1 }{ (2n+1)(2n+2)}=\log2,
$$
$$
\sum_{n=0}^\infty \frac{1 }{ (n+1)(2n+1)(4n+1)}=\frac{\pi }{ 3} \virgule
$$
et
 $$
\sum_{n=0}^\infty\frac {1 }{ (3n+1)(3n+2)(3n+3) } = \frac{\pi\sqrt{3}}{12} 
-\frac{1}{4}\log 3.
$$
De faon gŽnŽrale, quand  la fraction rationnelle $Q$ au dŽnominateur a uniquement des p™les simples et rationnels, la somme de la sŽrie est une  combinaison linŽaire   de logarithmes de nombres algŽ\-briques. En effet, d'aprs \cite{AdhikariSaradhaShoreyTijdemanIM-01} lemme 5, si $k_j$ et $r_j$ sont des entiers positifs avec $r_j\le k_j$ et si $c_j$ sont des nombres complexes, si la sŽrie
$$
S=\sum_{n=0}^\infty
\sum_{j=1}^m \frac{c_j}{k_jn+r_j}
$$
converge, alors sa somme est
$$
S=\sum_{j=1}^m  \frac{c_j}{k_j} \sum_{t=1}^{k_j-1} 
(1-\zeta_j^{-r_j t})\log (1-\zeta_j^{t}),
$$
$\zeta_j$ dŽsignant une racine primitive $k_j$-ime de l'unitŽ. Ainsi quand les nombres $c_j$ sont algŽbriques et que ce nombre $S$ est non nul, alors non seulement il est transcendant (d'aprs le thŽorme \ref{T:Baker}), mais en plus on en conna"t de bonnes mesures de transcendance (voir par exemple \cite{miwGL326-00} et \cite{AdhikariSaradhaShoreyTijdemanIM-01}). En particulier il n'est pas un nombre de Liouville.

D'autres exemples de sŽries de la forme  (\ref{E:serie}) prenant une valeur transcendante sont
$$
\sum_{n=1}^\infty \frac{1 }{ n^2}=\frac{\pi^2 }{ 6}
$$ 
et \cite{BundschuhMM-79}
\begin{equation}
\label{E:Bundschuh}
\sum_{n=0}^\infty\frac {1 }{ n^2+1}=\frac{1 }{ 2}
+\frac{\pi }{ 2}\cdot \frac{ e^\pi+e^{-\pi}
 }{ e^\pi-e^{-\pi}}\cdotp
\end{equation}
La transcendance de ce dernier nombre provient du thŽorme \ref{T:Nesterenko} de Nesterenko.

 Ces exemples soulvent plusieurs questions. Voici la premire
 
\begin{question} 
Quelles sont les pŽriodes parmi les nombres {\rm (\ref{E:serie})}?
   \end{question}

 Il y a beaucoup de pŽriodes parmi ces nombres   (\ref{E:serie}) (voir par exemple le lemme \ref{L:zeta} ci-dessous),
mais on s'attend plut™t ˆ ce qu'un nombre tel que (\ref{E:Bundschuh}) n'en soit pas une.

Sur la nature arithmŽtique des sŽries  (\ref{E:serie}), on peut espŽrer l'ŽnoncŽ suivant:
 
\begin{conjecture}
\label{C:Transcendancedeseries} 
Un nombre de la forme
 $$
  \sum_{n\ge 0 \atop   Q(n)\not=0}\frac {P(n) }{ Q(n)}
  $$
est soit rationnel, soit transcendant.
Ce n'est jamais un nombre de Liouville.
De plus, s'il est rationnel, alors la sŽrie est \og tŽlescopique\fg.

\end{conjecture}

Par \og sŽrie tŽlescopique\fg \ nous entendons une sŽrie dont la somme est rationnelle, la dŽmonstration de ce fait reposant sur  l'argument du lemme \ref{L:telescopique}.

\par

Le cas particulier des fractions rationnelles de la forme 
$$
P(X)/Q(X)=X^{-s}
$$ 
mŽrite une section spŽciale.

\section{Valeurs de la fonction zta de Riemann}

Commenons par un rŽsultat ŽlŽmentaire concernant les valeurs de la fonction zta  de Riemann aux entiers positifs.

\begin{lemme}
\label{L:zeta}
Pour $s\ge 2$
 $$
 \zeta(s)=\sum_{n\ge 1}\frac {1 }{ n^s}
 $$
est une pŽriode.
\end{lemme}

\noindent{\it DŽmonstration:} on vŽrifie facilement l'ŽgalitŽ
\begin{equation}
\label{E:MZVperiodes}
 \displaystyle
 \zeta(s)=\int_{1>t_1>\cdots>t_s>0}
 \frac{dt_1}{t_1}\cdots  \frac{dt_{s-1}}{t_{s-1}} \cdot \frac{dt_s}{1-t_s}
 \cdotp
\end{equation}

Il sera commode d'utiliser la notation des intŽgrales itŽrŽes de Chen (\cite{CartierB-00} \S~2.6) et d'Žcrire la relation (\ref{E:MZVperiodes}) sous la forme
\begin{equation}
\label{E:Chen}
  \zeta(s)=\int_0^1 \omega_0^{s-1}\omega_1
  \quad\hbox{avec}\quad
 \omega_0=\frac{dt }{ t}   \quad\hbox{et}\quad
 \omega_1=\frac{dt }{ 1-t} \cdotp
\end{equation}

La nature arithmŽtique des valeurs de la fonction zta de Riemann en des entiers positifs {\it pairs} est connue depuis  Euler:
$$
 \pi^{-2k}\zeta(2k)\in\BbbQ\quad\hbox{pour}\quad k\ge 1.
 $$
Ces nombres rationnels s'expriment en termes des nombres de Bernoulli \cite{CartierB-00}, formule (38). ƒtablir un rŽsultat de rationalitŽ (ou d'algŽ\-bricitŽ) de certains nombres est en gŽnŽral plus fŽcond que d'Žtablir des ŽnoncŽs d'irrationalitŽ ou de transcendance -- cependant notre propos est de faire le point sur les rŽsultats de transcendance: ce sont eux qui assurent que toute la richesse potentielle que reclent des relations algŽbriques entre les nombres considŽrŽs a bien ŽtŽ exploitŽe.

La principale question diophantienne que posent les nombres d'Euler est de savoir quelles relations algŽbriques existent entre les nombres 
$$
\zeta(2),\quad \zeta(3), \quad \zeta(5), \quad \zeta(7)\dots\ ?
$$
 
On conjecture qu'il n'y en a pas (\cite{CartierB-00}  et  \cite{FischlerB-02} Conjecture 0.1). Autrement dit

\begin{conjecture}\label{C:zeta} 
Les nombres
$$
\zeta(2),\quad \zeta(3), \quad \zeta(5), \quad \zeta(7)\dots\ 
$$
sont algŽbriquement indŽpendants.
\end{conjecture}

On sait trs peu de choses dans cette direction: le thŽorme de  Lindemann: affirme que le nombre $\pi$ est transcendant, donc aussi $\zeta(2k)$  pour tout  entier $k\ge 1$. En 1978 R.~ApŽry a dŽmontrŽ que le nombre   $\zeta(3)$ est irrationnel.  La dŽmonstration d'ApŽry permet de montrer que le nombre $\zeta(3)$ n'est pas un nombre de Liouville, la meilleure mesure d'irrationalitŽ Žtant celle de Rhin et Viola \cite{RhinViolaAA-01}:
$$
\left|\zeta(3)-\frac{p}{q}\right|>q^{-\mu} 
$$
pour $q$ suffisamment grand, avec $\mu=5,513\dots$

Les travaux rŽcents de T.~Rivoal, puis de  K.~Ball et  W.~Zudilin notamment, apportent les premires informations sur la nature arithmŽtique des valeurs de la fonction zta aux entiers impairs: par exemple l'espace vectoriel sur le corps des nombres rationnels engendrŽ par les nombres $\zeta(2k+1)$, $k\ge 1$ a une dimension infinie   (cf. \cite{FischlerB-02}).
 
Une Žtape prŽliminaire en vue d'une dŽmonstration de la conjecture  \ref{C:zeta} consiste ˆ linŽariser le problme: les mŽthodes diophantiennes sont en effet plus performantes pour Žtablir des ŽnoncŽs d'indŽpendance linŽaire (comme le thŽorme  \ref{T:Baker} de Baker, ou mme le thŽorme de Lindemann-Weiserstrass, qui peut s'Žnoncer de manire Žquivalente comme un rŽsultat d'indŽpendance algŽbrique (\cite{FeldmanNesterenko-98}  Th.~2.3') ou linŽaire (\cite{FeldmanNesterenko-98}  Th.~2.3) que pour Žtablir des rŽsultats d'indŽpendance algŽ\-brique. 
 
 Euler avait dŽjˆ remarquŽ que le produit de deux valeurs de la fonction zta (de Riemann comme on l'appelle maintenant!) Žtait  encore la somme d'une sŽrie.
En effet, 
de la relation
$$
\sum_{n_1\ge 1} n_1^{-s_1} \sum_{n_2\ge 1} {n_2}^{-s_2}=
\sum_{n_1>n_2\ge 1} n_1^{-s_1}{n_2}^{-s_2}+
\sum_{n_2>n_1\ge 1} n_2^{-s_2}{n_1}^{-s_1}+
\sum_{n\ge 1} n^{-s_1-s_2}
$$ 
on dŽduit, pour $s_1\ge 2$ et $s_2\ge 2$,
$$
\zeta(s_1)\zeta(s_2)=\zeta(s_1,s_2)+\zeta(s_2,s_1)+
\zeta(s_1+s_2)
$$
avec
$$\zeta(s_1,s_2)=
\sum_{n_1>n_2\ge 1} n_1^{-s_1}{n_2}^{-s_2}.
$$

Pour $k$, $s_1,\ldots,s_k$  entiers positifs avec $s_1\ge 2$, on pose
$\us=(s_1,\ldots,s_k)$ et 
\begin{equation}
\label{E:serieMZV}
\zeta(\us)=\sum_{n_1>n_2>\cdots>n_k\ge 1}\frac {1 }{
n_1^{s_1}\cdots n_k^{s_k}}\cdotp
\end{equation}
Ces nombres sont appelŽs \og valeurs zta multiples\fg, ou encore  \og MZV\fg ({\it Multiple Zeta Values}).
 Pour $k=1$ on retrouve bien entendu les nombres d'Euler $\zeta(s)$.

\begin{remarque} 
Chacun des nombres $ \zeta(\us)$ est une pŽriode: en effet, avec la notation {\rm (\ref{E:Chen})}  des intŽgrales itŽrŽes de Chen,   on a (cf. \cite{CartierB-00} \S~2.6):
\begin{equation}
\label{E:integraleMZV}
  \zeta(\us)=\int_0^1 \omega_0^{s_1-1}\omega_1\cdots
  \omega_0^{s_k-1}\omega_1.
\end{equation}
\end{remarque}

Le produit de sŽries (\ref{E:serieMZV}) est une combinaison linŽaire de telles sŽries. Par consŽquent l'espace vectoriel (sur $\BbbQ$ ou sur $\Qbar$) engendrŽ par les $\zeta(\us)$ est aussi une algbre sur ce corps. De plus  le produit de deux intŽgrales (\ref{E:integraleMZV}) est aussi une combinaison linŽaire de telles intŽgrales. Par diffŽrence on obtient  des relations linŽaires non triviales ˆ coefficients rationnels entre les MZV.  On en obtient de nouvelles --- comme celle d'Euler $\zeta(2,1)=\zeta(3) $ --- en y ajoutant les relations que l'on obtient en 
rŽgularisant les sŽries et les intŽgrales divergentes \cite{CartierB-00}. 

  Une description exhaustive des relations linŽaires entre les MZV devrait thŽoriquement permettre de dŽcrire du mme coup toutes les relations algŽ\-briques entre ces nombres, et en particulier de  rŽsoudre le problme \ref{C:zeta} de l'indŽpendance algŽbrique sur le corps $\BbbQ(\pi)$ des valeurs de la fonction zta aux entiers impairs.
  Le but est donc de dŽcrire toutes les relations linŽaires ˆ coefficients rationnels entre les MZV. 
Soit
$\Zgoth_p$ le $\BbbQ$-sous-espace vectoriel de $\BbbR$ engendrŽ par les nombres $\zeta(\us)$ pour $\us$ de \og poids\fg  $s_1+\cdots+s_k=p$, avec $\Zgoth_0=\BbbQ$ et
$\Zgoth_1=\{0\}$.

Voici la conjecture de  Zagier (conjecture (108) de \cite{CartierB-00})
sur la dimension $d_p$ de $\Zgoth_p$.

\begin{conjecture}
  \label{C:Zagier}
Pour  $p\ge 3$ on a 
$
d_p=d_{p-2}+d_{p-3}$: 
$$
(d_0,\ d_1,\ d_2,\ldots)=
(1,0,1,1,1,2,2,\ldots).
$$
  \end{conjecture}

Cette conjecture s'Žcrit aussi
$$
\sum_{p\ge 0} d_pX^p=\frac{1 }{ 1-X^2-X^3}\cdotp
$$
Un candidat pour tre une base de l'espace $\Zgoth_p$ est proposŽ par M.~Hoffman (\cite{HoffmanJA-97}, Conjecture C):

  \begin{conjecture}
  \label{C:Hoffman}
Une base de  $\Zgoth_p$ sur $\BbbQ$ est donnŽe par les nombres
$\zeta(s_1,\ldots,s_k)$, $s_1+\cdots+s_k=p$,  o chacun des $s_i$ est soit
$2$, soit
$3$.

  \end{conjecture}

Cette conjecture est compatible avec ce qui est connu  pour  $p\le 16$ (travaux de Hoang Ngoc
Minh notamment). Par exemple, en notant $\{a\}_b$ la suite formŽe de $b$ occurences de $a$,  les $7$ valeurs suivantes devraient tre une base de l'espace vectoriel 
$\Zgoth_{10}$:
\begin{eqnarray}   
&\zeta\bigl(\{2\}_5\bigr),\;
\zeta\bigl(\{2\}_2,\{3\}_2\bigr),\;
\zeta\bigl(\{2,3\}_2\bigr),\;
\zeta\bigl((2,\{3\}_2,2\bigr),\;
\nonumber\\
&\zeta\bigl(3,\{2\}_2,3\bigr),\;
\zeta\bigl(\{3,2\}_2\bigr),\;
\zeta\bigl(\{3\}_2,\{2\}_2\bigr).
\nonumber
\end{eqnarray}

\begin{exemple}
Voici les petites valeurs de $d_p$:
 
 \smallskip
 \noindent
 $\bullet$
$d_0=1$ car par convention $ \zeta(s_1,\ldots,s_k)=1$ 
pour $k=0$.

 \smallskip
 \noindent
 $\bullet$
$d_1=0$  car $\{(s_1,\ldots,s_k)\; ; \; 
s_1+\cdots+s_k=1,\; s_1\ge 2\}=\emptyset$.

 \smallskip
 \noindent
 $\bullet$
$d_2=1$  car  $\zeta(2)\not=0$

 \smallskip
 \noindent
 $\bullet$
$d_3=1$  car  $\zeta(2,1)=\zeta(3)\not=0$

 \smallskip
 \noindent
 $\bullet$
$d_4=1$  car  $\zeta(4)\not=0$
et
$$
\zeta(3,1)=\frac{1}{4}\zeta(4),\quad
\zeta(2,2)=\frac{3}{4}  \zeta(4),\quad
\zeta(2,1,1)= \zeta(4)=\frac{2}{5}\zeta(2)^2.
$$

\end{exemple}

La premire valeur de $d_p$ qui ne soit pas connue est $d_5$. La conjecture \ref{C:Zagier} donne  $d_5=2$, et on sait $d_5\in\{1,2\}$ car
\begin{eqnarray}   
 \zeta(2,1,1,1)&=&\displaystyle\zeta(5),\nonumber\\
\zeta(3,1,1)&=&\displaystyle\zeta(4,1)=\displaystyle2\zeta(5)-\zeta(2)\zeta(3),\nonumber\\
\zeta(2,1,2)&=&\displaystyle
 \zeta(2,3)=\displaystyle\frac{9 }{ 2}\zeta(5)-2\zeta(2)\zeta(3),\nonumber\\
\zeta(2,2,1)&=&\displaystyle\zeta(3,2)=\displaystyle3\zeta(2)\zeta(3)-\frac{11 }{ 2}\zeta(5).\nonumber
\end{eqnarray}

Donc  $d_5=2$  si et seulement si le nombre
$\zeta(2)\zeta(3)/\zeta(5)$ est irrationnel.

La conjecture \ref{C:Zagier} prŽdit une valeur exacte pour la dimension $d_p$ de $\Zgoth_p$. La question diophantienne est d'Žtablir la minoration. La majoration a ŽtŽ Žtablie rŽcemment gr‰ce aux travaux de A.B.~Goncharov   \cite{GoncharovECM-01} et
T.~Terasoma \cite{TerasomaIM-02} (voir aussi le thŽorme 6.4 de \cite{HoffmanJA-97}):

 {\it  Les entiers $\delta_p$ dŽfinis par la relation de rŽcurrence de la conjecture de Zagier  
$$
\delta_p=\delta_{p-2}+\delta_{p-3}
$$
avec les conditions initiales
$\delta_0=1$, $\delta_1=0$ fournissent une majoration pour la dimension $d_p$ de 
 $\Zgoth_p$. }

\section{Fonctions hypergŽomŽtriques}

Pour $a$, $b$, $c$ et $z$ nombres complexes avec $c\not\in\BbbZ_{\le 0}$ et  $|z|<1$, on dŽfinit la {\it fonction hypergŽomŽtrique de Gauss  } (voir par exemple \cite{FeldmanNesterenko-98}, Chap.~1 \S~3.6, Chap.~2 \S~3.2)
$$
{}_2 F_1\left({a},\
{b}\ ;\
{c}\ \bigl|\
{z}\right)=\sum_{n=0}^\infty \frac{(a)_n(b)_n }{ (c)_n}\cdot \frac{z^n }{ n!}
$$
o
$$
(a)_n=a(a+1)\cdots(a+n-1).
$$

\begin{exemple}
Si on note 
$K(z)$  l'intŽgrale elliptique de Jacobi de premire espce
$$
K(z)=\int_0^1 \frac{dx}{\sqrt{(1-x^2)(1-z^2x^2)}}\virgule
$$
$P_n$ le $n$-ime polyn™me de Legendre et 
$T_n$ le $n$-ime polyn™me de Chebyshev:
$$
P_n(z)=\frac{1}{n!}\left(\frac{d}{dz}\right)^n(1-z^2)^n,
\qquad
T_n(\cos z)=\cos (nz)
$$
 on a

\begin{equation}
\label{E:Hyper}
\left\{\begin{array}{cccccc}
{}_2 F_1\left({a},\ {1}\ ;\ {1}\ \bigl|\
{z}\right)\ &=& \ \displaystyle \
\frac{1}{(1-z_{\mathstrut})^a} \virgule
\\
{}_2 F_1\left(1,\ {1}\ ;\ {2}\ \bigl|\
{z}\right)\ &=& \ \displaystyle \
\frac{1}{z_{\mathstrut}}\log(1+z) \virgule
\\ 
{}_2 F_1\left({1/2},\ {1}\ ;\ {3/2}\ \bigl|\
{z^2}\right)\ &=& \ \displaystyle \
\frac{1}{2z_{\mathstrut}} \log \frac{1+z}{1-z} \virgule
\\
{}_2 F_1\left({1/2},\ {1/2}\ ;\ {3/2}\ \bigl|\
{z^2}\right)\ &=& \ \displaystyle \
\frac{1}{z_{\mathstrut}} \arcsin z,
\\
{}_2 F_1\left({1/2},\ {1/2}\ ;\ {1}\ \bigl|\
{z^2}\right)\ &=& \ \displaystyle \
\frac{2}{\pi_{\mathstrut}}K(z),
\\
{}_2 F_1\left({-n},\ {n+1}\ ;\ {1}\ \bigl|\
{(1+z)/2}\right)\ &=& \ \displaystyle \
2^{-n}P_n(z),
\\
{}_2 F_1\left({-n},\ {n}\ ;\ {1/2}\ \bigl|\
{(1+z)/2}\right)\ &=& \ \displaystyle \
(-1)^{n}T_n(z).
\end{array} \right.
\end{equation} 

\end{exemple}

Pour $c>b>0$ nombres rationnels, on a (Euler, 1748)
$$
{}_2 F_1\left({a},\
{b}\ ;\
{c}\ \bigl|\
{z}\right)=
\frac{\Gamma(c)}
{\Gamma(b)\Gamma(c-b)}
\int_0^1 t^{b-1} (1-t)^{c-b-1}
 (1-tz)^{-a}dt;
 $$
 de la formule de rŽflexion de la fonction Gamma, jointe ˆ la relation (\ref{E:Beta}) liant les fonctions Bta et Gamma, on dŽduit
 $$
 \frac{\Gamma(c)}
{\Gamma(b)\Gamma(c-b)}
\in\frac{1}{\pi}\calP.
$$
Il en rŽsulte    \cite{KontsevichZagierMU-00} \S~2.2 que pour $a$, $b$, $c$ rationnels avec $c\not\in\BbbZ_{\le 0}$ et $z\in\Qbar$ avec $|z|<1$,
\begin{equation}
\label{E:Fdeuxun}
{}_2 F_1\left({a},\
{b}\ ;\
{c}\ \bigl|\
{z}\right)
\in\frac{1}{\pi}\calP.
\end{equation}
Rappelons que
$\calP\subset (1/\pi)\calP$.
Il est suggŽrŽ dans  \cite{KontsevichZagierMU-00}  \S~2.2 que sous les mmes conditions, ${}_2 F_1\left({a},{b}\; ;{c} \bigl| {z}\right) $
n'appartient pas ˆ $\calP$.
 
\begin{remarque} 
Pour $a$, $b$, $c$ rŽels avec $c>a+b$ et $c\not\in\BbbZ_{\le 0}$, on a (Gauss)
 $$
 {}_2 F_1\left({a},\
{b}\ ;\
{c}\ \bigl|\
{1}\right)=
\frac{\Gamma(c-a)\Gamma(c-b)}
{\Gamma(c)\Gamma(c-a-b)}\cdotp
$$
\end{remarque}

\begin{remarque} 
Une relation remarquable faisant intervenir l'invariant modulaire $j(z)=J(e^{2i\pi z})$ et la sŽrie d'Eisenstein $E_4=Q$ (cf {\rm (\ref{E:Ramanujan})}) est la suivante {\rm
 \cite{KontsevichZagierMU-00}  \S~2.3}, due ˆ Fricke et Klein:
$$
{}_2 F_1\left(\frac{1 }{ 12}\virgule\
\frac{5 }{ 12}\ ;\
1\left|\
\frac{1728 }{ j(z)}\right)\right.=Q(z)^{1/4}.
$$

\end{remarque}

La transcendance des valeurs  ${}_2 F_1\left({a},
{b} ;{c} \bigl|{z}\right)$  des fonctions hypergŽomŽ\-tri\-ques quand $a$, $b$, $c$ et $z$ sont rationnels a ŽtŽ ŽtudiŽe ds 1929 par C.L.~Siegel  \cite{SiegelAMS-49}. On doit ˆ A.B.~Shidlovskii et ˆ son Žcole de nombreux rŽsultats sur la question (voir \cite{FeldmanNesterenko-98}).
 
En 1988, J.~Wolfart \cite{WolfartIM-88} a ŽtudiŽ  l'ensemble  $\calE$ des nombres algŽbriques $\xi$ tels que  ${}_2 F_1(\xi)$ soit aussi algŽbrique. Quand ${}_2 F_1$ est une fonction algŽbrique, $\calE=\overline {\BbbQ}$ est l'ensemble de tous les nombres algŽbriques. Supposons maintenant que ${}_2 F_1$ est une fonction transcendante.  Wolfart  \cite{WolfartIM-88} montre que l'ensemble $\calE$ est en bijection avec un ensemble de variŽtŽs abŽliennes de type CM -- il s'agit donc d'une extension en dimension supŽrieure du thŽorme de Th.~Schneider sur la transcendance de l'invariant modulaire $j$ (\cite{SchneiderMA-37}, \cite{Schneider-57} Th.~17).  

La dŽmonstration de Wolfart utilise le fait que les nombres   ${}_2 F_1(a,b,c;z)$ sont reliŽs aux pŽriodes de formes diffŽrentielles sur la courbe   
$$
y^N=x^A(1-x)^B(1-zx)^C
$$
avec  $A=(1-b)N$, $B=(b+1-c)N$, $C=aN$, tandis que $N$ est le plus petit dŽnominateur commun de $a,b,c$ (voir ˆ ce sujet  \cite{KoblitzRohrlichCJM-78}). L'outil transcendant est le   thŽorme du sous-groupe analytique de WŸstholz \cite{WustholzAnnM-89}.  Ë l'occasion de ces recherches,  F.~Beukers et J.~Wolfart \cite{BeukersWolfartD-86} ont Žtabli de nouvelles relations qui n'avaient pas ŽtŽ observŽes avant, comme 
$$
{}_2 F_1\left(\frac{1 }{ 12}\virgule\
\frac{5 }{ 12}\ ;\
\frac{1 }{ 2}\left|\
\frac{1323 }{ 1331}\right)\right.=\frac{3 }{ 4}\root {4} \of {11}
$$
et
$$
{}_2 F_1\left(\frac{1 }{ 12}\virgule\
\frac{7 }{ 12}\ ;\
\frac{2 }{ 3}\left|\
\frac{64000 }{64009}\right)\right.=\frac{2 }{ 3}\root {6} \of {253}.
$$

Quand le groupe de monodromie de l'Žquation diffŽrentielle hypergŽo\-mŽ\-trique satisfaite par ${}_2 F_1$ est un groupe triangulaire arithmŽtique,  l'ensemble   $\calE$ est infini. Wolfart  \cite{WolfartIM-88} conjectura que rŽciproquement,  l'ensemble   $\calE$ est fini si le groupe de monodromie n'est pas arithmŽtique. Les travaux de P.~Cohen et J.~Wolfart \cite{CohenWolfart-90}, puis de P.~Cohen et J.~WŸstholz  \cite{CohenWustholz-99} ont Žtabli un lien entre cette question et la conjecture d'AndrŽ-Oort \cite{AndreLivre-89,Oort-94}, selon laquelle les sous-variŽtŽs spŽciales de variŽtŽs de Shimura sont prŽcisŽment les sous-variŽtŽs qui contiennent un sous-ensemble Zariski dense de points spŽciaux. P.~Cohen a remarquŽ qu'un cas particulier de la conjecture d'AndrŽ-Oort en dimension $1$ suffit; le rŽsultat crucial a ŽtŽ Žtabli par B.~Edixhoven et A.~Yafaev \cite{EdixhovenYafaev-03}: {\it dans une variŽtŽ de Shimura, une courbe contient une infinitŽ de points appartenant ˆ la mme orbite de Hecke d'un  point spŽcial si et seulement si elle est de type Hodge}. Cela permet de rŽpondre ˆ la question initiale de C.L.~Siegel; de manire prŽcise,  en regroupant les rŽsultats de \cite{CohenWolfart-90, CohenWustholz-99, EdixhovenYafaev-03, WolfartIM-88}, on dŽduit: {\it l'ensemble exceptionnel $\calE$ est en bijection avec un ensemble de points dans la mme orbite de Hecke d'un point spŽcial (CM) sur une courbe dans une variŽtŽ de Shimura dŽfinie sur $\Qbar$.  L'ensemble   $\calE$ est infini si le groupe de monodromie est arithmŽtique, il  est fini si le groupe de monodromie n'est pas arithmŽtique.}

Ë la fonction hypergŽomŽtrique de Gauss on associe la fraction continue de Gauss
\begin{eqnarray}
G(z)&=&G(a,b,c;z) =
{}_2 F_1
\left({a},\
{b+1}\ ;\
{c+1}\ \bigl|\
{z}\right) /{}_2 F_1
\left({a},\
{b}\ ;\
{c}\ \bigl|\
{z}\right) 
\nonumber\\
&=&1/\Bigl(1-g_1z/\bigl(1-g_2z/(\cdots)\bigr)\Bigr)
\nonumber
\end{eqnarray}
ˆ coefficients
\begin{eqnarray}
g_{2n-1}&=&(a+n-1)(c-b+n-1)/\bigl((c+2n-2)(c+2n-1)\bigr),
\nonumber\\
g_{2n}&=& (b+n)(c-a+n)/\bigl((c+2n-1)(c+2n)\bigr).
\nonumber
\end{eqnarray}
J.~Wolfart \cite{WolfartBudapest-90} a montrŽ que si les paramtres  $a,b,c$ sont rationnels, $c\not=0$, $-1,-2,\cdots$, et si $G(z)$ n'est pas une fonction algŽbrique, alors pour presque toutes les valeurs algŽbriques de l'argument $z$  la valeur de  $G(z)$ est transcendante. Il utilise le thŽorme de G. WŸstholz \cite{WustholzAnnM-89} ainsi que des rŽsultats de   G. Shimura et Y. Taniyama sur les variŽtŽs abŽliennes. 

Pour les nombres dont nous venons de parler, liŽs aux fonctions hypergŽo\-mŽtriques, dont on sait Žtablir la transcendance, on conna"t Žgalement des mesures d'approximation par des nombres algŽbriques (voir \cite{FeldmanNesterenko-98} et \cite{GaudronCRAS-01, GaudronFLLEVA-04} notamment).

Les fonctions hypergŽomŽtriques de Gauss sont des cas particuliers d'une famille plus vaste, formŽe des {\it  fonctions hypergŽomŽtriques gŽnŽralisŽes} (voir par exemple \cite{FeldmanNesterenko-98} Chap.~2, \S~6).

Pour $p$ entier $\ge 2$ et $a_1,\ldots,a_p$, $b_1,\ldots,b_{p-1}$  et $z$ nombres complexes avec $b_i\not\in\BbbZ_{\le 0}$ et  $|z|<1$, on dŽfinit 
$$
{}_p F_{p-1}\left(
{{a_1\ ,\ \ldots\ ,\ a_p}\atop {b_1,\ldots,b_{p-1}}}
\ \Bigl|\
{z}\right)=\sum_{n=0}^\infty \frac{(a_1)_n\cdots (a_p)_n }{
(b_1)_n\cdots (b_{p-1})_n}\cdot \frac{z^n }{ n!}\cdotp
$$

\begin{exemple} Les fonctions
$$
{}_1 F_{0}\left(
{{1/n }\atop { }}
\ \Bigl|\ z^n\right)=
\root n \of {1-z^n}
$$
et
$$
{}_3 F_{2}\left(
{{1/4, 1/2, 3/4}\atop {1/3,2/3}}
\ \Bigl|\
\frac{2^8z}{3^3}\right)=
\sum_{k=0}^\infty {{4k}\choose {k}}z^k
$$
sont algŽbriques. La fonction de Bessel
$$
J_0(z)=
{}_0 F_{1}\left(
{{ }\atop {1 }}
\ \Bigl|\ \frac{-z^2}{4}\right) =\sum_{n\ge 0}(-1)^n\frac{z^{2n}}{2^{2n}(n!)^2}
$$
est transcendante.
\end{exemple}

Par rŽcurrence sur $p$ ˆ partir de (\ref{E:Fdeuxun}) 
 on dŽmontre:

\begin{propos}
\label{P:Hypergeometrique}
Pour  $a_1,\ldots,a_p$, $b_1,\ldots,b_{p-1}$  nombres rationnels avec  $p\ge 2$, $b_i\not\in\BbbZ_{\le 0}$ et pour $z\in\Qbar$ avec $|z|<1$, on a
$$
{}_p F_{p-1}\left(
{{a_1\ ,\ \ldots\ ,\ a_p}\atop {b_1,\ldots,b_{p-1}}}\ 
\Bigl|\
{z}\right)
\in\frac{1}{\pi^{p-1}}\calP.
$$
\end{propos}

On peut consulter \cite{FeldmanNesterenko-98} pour conna"tre l'Žtat de la question de la nature arithmŽtique des valeurs de fonctions hypergŽomŽtriques gŽnŽralisŽes, aussi bien sous l'aspect qualitatif (irrationalitŽ, transcendance, indŽpendance algŽ\-brique) que quantitatif (mesures d'approximation, de transcendance ou d'indŽ\-pendance algŽbrique).

\section{Mesure de Mahler de polyn™mes en plusieurs variables}

Soit $P\in\BbbC[z_1,\ldots,z_n,z_1^{-1},\ldots,z_n^{-1}]$  un polyn™me de Laurent non nul en $n$ variables. On dŽfinit la mesure de Mahler $\rmM(P)$ et la mesure de Mahler logarithmique $\mu(P)$ par
$$
\mu(P)= \log \rmM(P)=
\int_0^1 \cdots\int_0^1 \log
\left| P(e^{2i\pi t_1},\ldots,e^{2i\pi t_n} ) \right| dt_1\cdots dt_n.
$$

Dans le cas le plus simple $n=1$ on Žcrit
$$
P(z)=\sum_{i=0}^d a_{d-i}z^i=a_0\prod_{i=1}^d (z-\alpha_i)
$$
et on a
$$
\rmM(P)=|a_0|\prod_{i=1}^d \max\{1,|\alpha_i|\}.
$$
On en dŽduit notamment $\mu(P)\in\calP$.

Plus gŽnŽralement, 
pour $P\in\Qbar[z_1,\ldots,z_n,z_1^{-1},\ldots,z_n^{-1}]$ on a
$$
\mu(P)\in\frac{1}{\pi^n}\calP.
$$
D.~Boyd et C.J.~Smyth (voir les rŽfŽrences dans \cite{BoydZakopane-99}) ont calculŽ un certain nombre d'exemples de valeurs de la fonction $\mu$ qu'ils ont exprimŽes en termes de valeurs spŽciales de fonctions $L$ attachŽes ˆ des caractres de Dirichlet. D.~Boyd et F.~Rodriguez Villegas ont obtenu des rŽsultats du mme genre faisant intervenir des fonctions $L$ attachŽes ˆ des courbes elliptiques. Ensuite D.~Boyd, F.~Rodriguez Villegas, V.~Maillot et S.~Vandervelde ont exprimŽ certaines mesures de Mahler logarithmiques ˆ l'aide de combinaisons de valeurs de la fonction dilogarithme. Mis ˆ part le cŽlbre exemple
$$
\mu(1+z_1+z_2+z_3)=\frac{7}{2\pi^2}\zeta(3).
$$ 
dž ˆ C.J.~Smyth, les rŽsultats connus concernent principalement le cas de deux variables; cependant les travaux de C.~Deninger donnent des espoirs pour le cas gŽnŽral.

\section{Exponentielles de pŽriodes et pŽriodes exponentielles}

\subsection{Exponentielles de pŽriodes}

Parmi les suggestions de Kontsevich et Zagier dans \cite{KontsevichZagierMU-00}, on relve dans le \S~1.2 celle qui prŽdit que les nombres $1/\pi$ et $e$ ne sont pas des pŽriodes. Parmi les candidats ˆ ne pas tre des pŽriodes on peut ajouter 
$ e^\pi $ et $ e^{\pi^2}$.
 
On conna"t la transcendance de $e^\pi$ (A.O.~Gel'fond, 1929 -- c'est une 
 consŽ\-quence du thŽorme \ref{T:Baker}), mais pas celle de $e^{\pi^2}$.

  \begin{conjecture}
  \label{C:troislogs} 
 Soient $\alpha_1,\alpha_2,\alpha_3$ des nombres algŽbriques non nuls. Pour $j=1,2,3$ soit $\log\alpha_j\in\BbbC\setminus\{0\}$ un logarithme non nul de $\alpha_j$, c'est-ˆ-dire un nombre complexe non nul tel que  $e^{\log\alpha_j}=\alpha_j$.  Alors
 $$
( \log\alpha_1)\log\alpha_2)\not=\log\alpha_3.
 $$
   \end{conjecture}

\begin{exemple}   
Avec $\log\alpha_1=\log\alpha_2=i\pi$ on dŽduit la
 transcendance du nombre $e^{\pi^2}$. Un autre exemple est la
transcendance du nombre  $2^{\log 2}$.
 \end{exemple}
 
 D'autres conjectures sont proposŽes dans \cite{miwHyderabad-03}, ˆ la fois pour la fonction exponentielle (conjectures des trois, quatre, cinq exponentielles) et pour les fonctions elliptiques. Des cas trs particuliers de ces conjectures sont Žtablis. L'appendice de \cite{miwHyderabad-03} par H.~Shiga Žtablit un lien avec des pŽriodes de surfaces de Kummer.

Les rŽsultats partiels que l'on conna"t sur la conjecture \ref{C:troislogs} et les questions autour de la conjecture des quatre exponentielles \cite{LangITN-66, miwGL326-00, miwHyderabad-03}   reposent sur la mŽthode de transcendance qui a permis ˆ Th. Schneider de rŽsoudre le septime problme de Hilbert en 1934. On peut noter que cette mŽthode fait jouer un r™le essentiel au thŽorme d'addition algŽbrique de la fonction exponentielle, ˆ savoir $e^{x+y}=e^xe^y$, et n'utilise pas l'Žquation diffŽrentielle de cette fonction, contrairement ˆ ce qui est suggŽrŽ ˆ la fin du \S~2.4 de \cite{KontsevichZagierMU-00}. Une autre mŽthode de transcendance qui ne fait pas intervenir de dŽrivations est celle de Mahler qui fait l'objet du livre de K.~Nishioka \cite{NishiokaLN}. Notons ˆ ce propos que la conjecture 5.4 de \cite{miwODPJMS-04} sur la  transcendance de nombres dont le dŽveloppement dans une base est donnŽ par une suite \og automatique\fg, qui faisait l'objet de travaux de J.H.~Loxton et A.J.~van der Poorten utilisant la mŽthode de Mahler, vient d'tre rŽsolue par B.~ Adamczewski,  Y.~Bugeaud  et F.~Lucas \cite{AdamczewskiBugeaudLucas-04, AdamczewskiBugeaud-04}, gr‰ce ˆ une mŽthode entirement diffŽrente de celle de Mahler, basŽe sur le thŽorme du sous-espace de W.M.~Schmidt. 

La conjecture \ref{C:troislogs} est un cas trs particulier de la conjecture selon laquelle {\it des logarithmes $\BbbQ$-linŽairement indŽpendants de nombres algŽbriques sont algŽbriquement indŽpendants.} Un des exemples les plus importants de nombres qui s'expriment comme valeur d'un polyn™me en des logarithmes de nombres algŽbriques est celui de dŽterminants de matrices dont les coefficients sont de tels logarithmes. Certains rŽgulateurs sont de cette forme, et dŽcider s'ils sont nuls ou non peut tre considŽrŽ comme un problme de transcendance. Quand ils ne sont pas nuls on conjecture que ces dŽterminants sont transcendants, et que ce ne sont pas de nombres de Liouville.

\subsection{PŽriodes exponentielles}

La dŽfinition suivante est donnŽe dans  \cite{KontsevichZagierMU-00} \S~4.3. Il est prŽcisŽ dans l'introduction de   \cite{KontsevichZagierMU-00} que la dernire partie de ce texte est l'{\oe}uvre uniquement du premier auteur.

\noindent{\bf DŽfinition}
 Une {\it pŽriode exponentielle} est une intŽgrale absolument convergente du produit d'une fonction algŽbrique avec l'exponentielle d'une fonction algŽbrique, sur un ensemble semi-algŽbrique, o tous les polyn™mes intervenant dans la dŽfinition ont des coefficients algŽbriques.

\begin{exemple} 
\noindent
     
  Dans l'algbre des pŽriodes exponentielles, on trouve Žvidemment les pŽ\-riodes, mais aussi les nombres 
  $$
  e^\beta=\int_{-\infty}^\beta e^xdx
  $$
  quand $\beta$ algŽbrique, le nombre
  $$
\sqrt{\pi}=\int_{-\infty}^\infty e^{-x^2} dx,
$$
les valeurs de la fonction Gamma aux points rationnels:
$$
\Gamma(s)=\int_0^\infty e^{-t}t^s \cdot \frac{dt }{ t}\virgule
$$
ainsi que  les valeurs des fonctions de Bessel aux points algŽbriques
$$
J_n(z)=
 \int_{|u|=1} 
\exp\left( \frac{z }{ 2}\Bigl(u-\frac{1 }{ u}\Bigr)\right)\frac {du }{ u^{n+1}
}\cdotp 
$$

\end{exemple}

Ces exemples sont interprŽtŽs par S.~Bloch et H.~Esnault {\rm \cite{BlochEsnault-99} } comme des pŽriodes issues d'une dualitŽ entre cycles homologiques et formes diffŽrentielles pour des connections ayant des points singuliers irrŽguliers sur des surfaces de Riemann.


\subsection{Constante d'Euler}

On ne sait pas dŽmontrer que le nombre  
$$
\gamma=\lim_{n\rightarrow\infty} \left(1+\frac{1 }{ 2}+\frac{1 }{ 3}+\cdots
+\frac{1 }{ n}-\log n\right)=0.5772157\dots
$$
  est irrationnel \cite{SondowPAMS-03}, mais on attend mieux:

  \begin{conjecture}
  \label{C:gammanonperiode}
  Le nombre $\gamma$ est transcendant.
    \end{conjecture}
    
    Un rŽsultat encore plus fort est suggŽrŽ par Kontsevich et Zagier \cite{KontsevichZagierMU-00} \S~1.1 et \S~4.3:

\begin{conjecture} 
  \label{C:Trdcegamma}
  Le nombre  $\gamma$ n'est pas une pŽriode, ni mme une pŽriode exponentielle.
  \end{conjecture}


\subsection{Un
analogue en dimension $2$ de la constante d'Euler}

  Pour chaque
$k\ge 2$,  dŽsignons par  $A_k$ l'aire minimale d'un disque fermŽ de  $\BbbR^2$ contenant au moins $k$ points de $\BbbZ^2$. Pour $n\ge 2$, posons \cite{GramainWeberMC-85}
$$
\delta_n=-\log n+\sum_{k=2}^n \frac{1 }{ A_k}\quad
  \hbox{et}\quad
  \delta=\lim_{n\rightarrow \infty}\delta_n.
 $$

F.~Gramain conjecture:
  
  \begin{conjecture}
  \label{C:Gramain}
$$
\delta=1+\frac{4 }{\pi}\bigl(\gamma L'(1)+ L(1)\bigr),
$$ 
  o  $\gamma$ est la constante d'Euler  et 
$$
L(s)=\sum_{n\ge 0}(-1)^n (2n+1)^{-s}.
$$
est la fonction $L$ du corps quadratique $\BbbQ(i)$ (fonction Bta de Dirichlet). 

\end{conjecture}
 
Comme
$
L(1)={\pi/ 4}
$ et 
\begin{eqnarray} 
\nonumber
\displaystyle
L'(1)&
\displaystyle
=&\displaystyle
\sum_{n\ge 0}(-1)^{n+1}\cdot\frac{\log(2n+1) }{  2n+1}\hfill\quad
\nonumber
\\ 
&
\displaystyle
=&\displaystyle
\frac{\pi }{ 4}\bigl(3\log
 \pi+2\log 2+\gamma-4\log\Gamma(1/4)\bigr),
 \nonumber
 \end{eqnarray}  
la conjecture \ref{C:Gramain} s'Žcrit aussi
$$
\delta=1+3\log
 \pi+2\log 2+2\gamma-4\log\Gamma(1/4)
=1.82282524\dots 
$$ 
Le meilleur encadrement connu pour $\delta$  est \cite{GramainWeberMC-85} 
$$
1.811\dots<\delta<1.897\dots
$$
Il semble vraisemblable que ce nombre $\delta$ n'est pas une pŽriode (et par consŽ\-quent est transcendant), mais, Žtant donnŽ le peu d'information que nous avons sur lui,  ce n'est probablement pas le meilleur candidat pour rŽsoudre la question de \cite{KontsevichZagierMU-00} (\S~1.2 problem 3) qui consiste ˆ exhiber un nombre qui n'est pas une pŽriode!


\section{CaractŽristique finie}

Les questions diophantiennes concernant les nombres complexes ont des analogues dans les corps de fonctions en caractŽris\-tique finie qui ont  fait l'objet de nombreux travaux \cite{Goss-92}. Les premiers outils dŽveloppŽs par L.~Carlitz (1935) ont ŽtŽ utilisŽs par I.I.~Wade (1941) qui a obtenu les premiers ŽnoncŽs de transcendance. Aprs divers travaux, notamment de J.M.~Geijsel et P.~Bundschuh en 1978, la thŽorie a ŽtŽ dŽveloppŽe de manire approfondie par Jing Yu ˆ partir des annŽes 1980, d'abord dans le cadre des modules elliptiques qui avaient ŽtŽ introduits par V.G.~Drinfel'd en 1974, ensuite dans le cadre des $t$-motifs de G.~Anderson ˆ partir de 1986. Pendant longtemps les rŽsultats en caractŽristique finie Žtaient des analogues des rŽsultats classiques relatifs aux nombres complexes, jusqu'ˆ ce que Jing Yu obtiennent des ŽnoncŽs qui vont plus loin que leurs analogues complexes \cite{JingYu-92}. 

L'utilisation, introduite dans ce contexte par L.~Denis en 1990, de la dŽrivation par rapport ˆ la variable du corps de fonctions, produit des ŽnoncŽs qui n'ont pas d'Žquivalents dans le cas classique des nombres complexes. Une autre particularitŽ de la caractŽristique finie est la possibilitŽ de considŽrer des produits tensoriels, permettant parfois de ramener des questions d'indŽpen\-dance algŽbrique ˆ des problmes d'indŽpendance linŽaire (un bel exemple est donnŽ par S.~David et L.~Denis dans \cite{DavidDenis-02}).

Deux exposŽs de synthse sur ce thme sont donnŽs, l'un en 1992  par Jing Yu dans \cite{JingYu-92}, l'autre en 1998 par W.D.~Brownawell \cite{BrownawellTrichy-98}.


\subsection{Un
analogue en caractŽristique finie de la constante d'Euler}

Un rŽsultat remarquable en caractŽristique finie est la transcendance de l'ana\-logue de la constante d'Euler. Le nombre complexe $\gamma$ peut tre dŽfini par
$$
\gamma=\lim_{s\rightarrow 1}\left(\zeta(s)-\frac{1 }{ s-1}\right)
$$
quand $\zeta$ dŽsigne la fonction zta de Riemann:
$$
{\zeta}(s)=\prod_{p} (1-p^{-s})^{-1}.
$$
Dans ce produit $p$ dŽcrit l'ensemble de nombres premiers. En caractŽristique finie le produit correspondant porte sur les polyn™mes irrŽductibles unitaires ˆ coefficients dans un corps fini $\BbbF_q$ ˆ $q$ ŽlŽments. Dans ce cas le produit converge au point $s=1$, et la valeur en ce point est donc un analogue de la constante d'Euler (dans un complŽtŽ $C$ d'une cl™ture algŽbrique de $\BbbF_q((1/T))$). Cet ŽlŽment de $C$ est transcendant sur $\BbbF_q(T)$: cela a ŽtŽ dŽmontrŽ par 
G.W.~Anderson et D.~Thakur en 1990 \cite{AndersonThakurAM-01}, mais ils remarquent que les outils dont disposait I.I.~Wade auraient suffit pour  Žtablir le rŽsultat ds 1940.


\subsection{La fonction Gamma de Thakur}

Par analogie avec le produit infini (\ref{E:GammaProduit}) dŽfinissant la fonction Gamma d'Euler, D.~Thakur dŽfinit (cf. \cite{BrownawellTrichy-98})
$$
\Gamma(z)=  z^{-1} \prod_{n\in A_+}^\infty \left(
1+\frac{z}{n}\right)^{-1}
$$
o $A=\BbbF_q[T]$ dŽsigne l'anneau des polyn™mes ˆ coefficients dans  $\BbbF_q$   et $A_+$ l'ensemble des polyn™mes unitaires. Cette fonction Gamma est mŽromorphe sur  $C$. Elle  satisfait des relations analogues aux relations standard satisfaites par la fonction Gamma d'Euler. De mme, le pendant en caractŽristique finie des relations de Deligne-Koblitz-Ogus a ŽtŽ Žtabli par Deligne, Anderson et Thakur (voir \cite{BrownawellPapanikolasCrelle-02}).  

En 1992 G. Anderson a introduit une notion de fonction soliton, qui, selon  \cite{BrownawellPapanikolasCrelle-02},  est un analogue en dimension supŽrieure de la fonction shtuka pour les modules de Drinfeld de rang $1$. Les fonctions mŽromorphes que  R.~Coleman avait utilisŽes pour Žtudier les endomorphismes de Frobenius sur les courbes de Fermat et d'Artin-Schreier avaient ŽtŽ interprŽtŽes par Thakur en termes de sa fonction Gamma. C'est en dŽveloppant le parallle avec certains ŽlŽments de la thŽorie des Žquations aux dŽrivŽes partielles que G. Anderson a ŽtŽ amenŽ ˆ introduire ses solitons. Cette thŽorie a ŽtŽ utilisŽe  en 1997 par S.K.~Sinha, qui a construit des $t$-modules ayant des pŽriodes dont les coordonnŽes sont des multiples par un nombre algŽbrique de valeurs de la fonction Gamma de Thakur en des points rationnels $a/f$, avec $a$ et $f$ dans $A_+$. Gr‰ce aux rŽsultats de transcendance de Jing Yu,  S.K.~Sinha a pu ainsi obtenir la transcendance de certaines valeurs de la fonction Gamma de Thakur en des points rationnels.

Ces travaux ont ŽtŽ poursuivis par W.D.~Brownawell et M.A.~Papanikolas qui montrent dans \cite{BrownawellPapanikolasCrelle-02} que les relations linŽaires ˆ coefficients algŽbriques entre les valeurs de la fonction Gamma de Thakur sont celles qui rŽsultent de Deligne-Anderson-Thakur. On peut considŽrer qu'il s'agit de l'analogue pour les corps de fonctions sur un corps fini du thŽorme de Wolfart et WŸstholz \cite{WolfartWustholzMA-85} sur l'indŽpendance linŽaire des valeurs de la fonction Bta. Ce qui est spŽcialement remarquable est qu'il est possible d'aller plus loin: 
dans \cite{AndersonBrownawellPapanikolasAM-04}, G.W.~Anderson, W.D.~Brownawell et M.A.~Papanikolas montrent  que {\it toutes les relations de dŽpendance algŽbrique  entre les valeurs de la fonction Gamma de Thakur rŽsultent des relations de Deligne-Anderson-Thakur.}

Pour les valeurs de la fonction Gamma la situation en caractŽristique finie est donc bien plus en avance que dans le cas classique o la conjecture de Rohrlich-Lang semble hors d'atteinte.

Tout rŽcemment M.~Papanikolas a Žtabli l'analogue pour les modules de Drinfeld de la conjecture sur l'indŽpendance algŽbrique de logarithmes de nombres algŽbriques (cf.~\S~8.1).

\section*{Remerciements}

Ce texte rŽsume des exposŽs donnŽs en 2003 et 2004 ˆ Sharhood (Iran), Taipei (Taiwan), Bangalore (Inde), Oulu (Finlande), SŽoul (CorŽe), Mahdia (Tunisie) et Besanon (France). Merci en particulier ˆ Khalifa Trimche pour son invitation au congrs annuel de la SociŽtŽ MathŽmatique de Tunisie en mars 2004 qui donne lieu ˆ cette rŽdaction.

\def\cprime{$'$} \def\cprime{$'$} \def\cprime{$'$}
\providecommand{\bysame}{\leavevmode ---\ }
\providecommand{\og}{``}
\providecommand{\fg}{''}
\providecommand{\smfandname}{\&}
\providecommand{\smfedsname}{\'eds.}
\providecommand{\smfedname}{\'ed.}
\providecommand{\smfmastersthesisname}{M\'emoire}
\providecommand{\smfphdthesisname}{Th\`ese}

\end{document}